\documentclass[11pt]{amsart}
\usepackage{amssymb, amsmath, amscd, eucal, rotating}

\input xy
 \xyoption{all}

\numberwithin{equation}{section}

\theoremstyle{plain}

\newtheorem{proposition}{Proposition}[section]
\newtheorem{theorem}[proposition]{Theorem}		
\newtheorem{corollary}[proposition]{Corollary}
\newtheorem{lemma}[proposition]{Lemma}

\theoremstyle{definition}

\newtheorem{definition}[proposition]{Definition}
\newtheorem{remark}[proposition]{Remark}

\newcommand{\C}{\mathbb C}
\newcommand{\R}{\mathbb R}
\newcommand{\ZBbb}{\mathbb Z}

\newcommand{\Gr}{\mathop{\rm Gr}\nolimits}

\newcommand{\codim}{\mathop{\rm codim}\nolimits}

\newcommand{\End}{\mathop{\rm End}\nolimits}
\newcommand{\coker}{\mathop{\rm coker}\nolimits}

\newcommand{\Hom}{\mathop{\rm Hom}\nolimits}

\newcommand{\GL}{\mathsf{GL}}

\newcommand{\U}{\mathsf{U}}

\newcommand{\SO}{\mathsf{SO}}

\DeclareMathOperator{\Rep}{Rep}
\DeclareMathOperator{\Ext}{Ext}
\DeclareMathOperator{\Vect}{Vect}
\DeclareMathOperator{\id}{id}
\DeclareMathOperator{\im}{im}
\DeclareMathOperator{\rank}{rank}
\DeclareMathOperator{\slope}{slope}

\begin{document}


\title[Homotopy groups of moduli spaces of stable quiver representations]
	{Homotopy groups of moduli spaces of stable quiver representations}

\author[Graeme Wilkin]{Graeme Wilkin}

\address{Department of Mathematics \\
		University of Colorado \\
		Boulder, CO 80303}

\thanks{}

\email{graeme.wilkin@colorado.edu}

\begin{abstract}
The purpose of this paper is to describe a method for computing homotopy groups of the space of $\alpha$-stable representations of a quiver with fixed dimension vector and stability parameter $\alpha$. The main result is that the homotopy groups of this space are trivial up to a certain dimension, which depends on the quiver, the choice of dimension vector, and the choice of parameter. As a corollary we also compute low-dimensional homotopy groups of the moduli space of $\alpha$-stable representations of the quiver with fixed dimension vector, and apply the theory to the space of non-degenerate polygons in three-dimensional Euclidean space.
\end{abstract}


\date{\today}



\maketitle

\thispagestyle{empty}


\baselineskip=16pt
\setcounter{footnote}{0}


\section{Introduction}

Symplectic quotients are important examples of symplectic manifolds, which arise in many areas of mathematics. When the quotient is a smooth manifold there exists a well-developed research program to compute the cohomology of these quotients via equivariant Morse theory, begun by Atiyah and Bott in \cite{AtiyahBott83}, and Kirwan in \cite{Kirwan84}. When the quotient is singular then in many cases of interest one would like to compute topological invariants of the smooth locus of this space, for which new techniques are needed. In the case of holomorphic bundles over a compact Riemann surface, an approach due to Daskalopoulos and Uhlenbeck is to use transversality theory to compute low-dimensional homotopy groups of the moduli space of stable holomorphic bundles. In this paper we extend the Daskalopoulos-Uhlenbeck techniques to the space of stable representations of a quiver, and compute homotopy groups of these moduli spaces.


Let $\Rep(Q, {\bf v})$ denote the space of representations of a quiver with fixed dimension vector ${\bf v}$, and let $G_{\bf v}$ be the associated group acting on $\Rep(Q, {\bf v})$ (see Section \ref{sec:background} for the precise definition). For a given stability parameter $\alpha \in \mathfrak{g}_{\bf v}^*$, denote the subset of $\alpha$-stable points by $\Rep(Q, {\bf v})^{\alpha-st}$, and the $\alpha$-semistable points by $\Rep(Q, {\bf v})^{\alpha-ss}$. In a typical situation, for almost all choices of $\alpha$ we have $\Rep(Q, {\bf v})^{\alpha-st} = \Rep(Q, {\bf v})^{\alpha-ss}$, and since the isotropy group of each $\alpha$-stable point is the subgroup consisting of scalar multiples of the identity, then there is a smooth structure on the quotient $\Rep(Q, {\bf v})^{\alpha-ss} / \! /_{\alpha} G_{\bf v}^\C = \Rep(Q, {\bf v})^{\alpha-st} / G_{\bf v}^\C$. 


Various methods have already been employed to compute topological invariants of the quotient $\Rep(Q, {\bf v})^{\alpha-st} / G_{\bf v}^\C$ in the case where the spaces $\Rep(Q, {\bf v})^{\alpha-ss}$ and $\Rep(Q, {\bf v})^{\alpha-st}$ co-incide. The major result in this direction is that of Reineke in \cite{Reineke03}, who used an argument involving counting rational points (analogous to the work of Harder and Narasimhan in \cite{HarderNarasimhan75} for the moduli space of stable vector bundles) to compute Betti numbers of the quotient. When $Q$ is the star-shaped polygon quiver studied in Section \ref{subsec:polygon}, Klyachko in \cite{Klyachko94} computed Betti numbers of $\Rep(Q, {\bf v})^{\alpha-st} / G_{\bf v}^\C$ for the case where $\Rep(Q, {\bf v})^{\alpha-ss} = \Rep(Q, {\bf v})^{\alpha-st}$, and Kapovich and Millson in \cite{KapovichMillson96} studied the geometry of $\Rep(Q, {\bf v})^{\alpha-st} / G_{\bf v}^\C$.

More recently, the results of \cite{haradawilkin} show that the Morse theory methods of Kirwan in \cite{Kirwan84} apply to the non-compact space $\Rep(Q, {\bf v})$, and as a consquence there is a surjective map in equivariant cohomology
\begin{equation*}
\kappa : H_{G_{\bf v}}^*\left(\Rep(Q, {\bf v}) \right) \twoheadrightarrow H_{G_{\bf v}}^*\left(\Rep(Q, {\bf v})^{\alpha-ss} \right) .
\end{equation*}
When the spaces $\Rep(Q, {\bf v})^{\alpha-ss}$ and $\Rep(Q, {\bf v})^{\alpha-st}$ co-incide then there is an additional isomorphism 
\begin{equation*}
H_{G_{\bf v}}^*(\Rep(Q, {\bf v})^{\alpha-ss}) = H_{G_{\bf v}}^*(\Rep(Q, {\bf v})^{\alpha-st}) \cong H^* (\mathcal{M}_\alpha(Q, {\bf v})^{st}) \otimes H^*(BU(1)) ,
\end{equation*}
and from the Morse theory one can recover the Betti number calculations of \cite{Reineke03}.


In general, when  the spaces $\Rep(Q, {\bf v})^{\alpha-ss}$ and $\Rep(Q, {\bf v})^{\alpha-st}$ do not co-incide then there are three different quotient varieties of interest: the GIT quotient $\mathcal{M}_\alpha (Q, {\bf v}) := \Rep(Q, {\bf v}) / \! /_\alpha G^\C$, the resolution of singularities $\tilde{\mathcal{M}}_\alpha (Q, {\bf v})$ of the GIT quotient (as studied by Kirwan in \cite{Kirwan85}), and the smooth locus $\mathcal{M}_\alpha(Q, {\bf v})^{st} := \Rep(Q, {\bf v})^{\alpha-st} / G^\C \subset \mathcal{M}_\alpha(Q, {\bf v})$. The object of study in this paper is the space $\mathcal{M}_\alpha(Q, {\bf v})^{st}$, and the results derived here on the topology of this space cannot be reproduced by the methods of \cite{Reineke03} or \cite{haradawilkin} when $\Rep(Q, {\bf v})^{\alpha-ss} \neq \Rep(Q, {\bf v})^{\alpha-st}$.


The moduli space of semistable holomorphic bundles of degree zero over a compact Riemann surface is an example of a singular K\"ahler quotient for which one would like to study the smooth locus: the moduli space of stable bundles. In the rank $2$ case Kirwan applied her general theory for computing intersection Betti numbers of GIT quotients to obtain results about the low-dimensional Betti numbers (see \cite{Kirwan86-2} and Remark 3.4 in \cite{Kirwan86}), and later Daskalopoulos in \cite{Daskal92} used infinite-dimensional methods in the spirit of Atiyah and Bott to obtain results about the homotopy groups of this space. Following this, Daskalopoulos and Uhlenbeck in \cite{DaskalUhlenbeck95} used an infinite-dimensional transversality argument to extend these results to the moduli space of stable bundles for any rank and degree, and obtained a simple formula for the homotopy groups of the moduli space of stable bundles up to a dimension bound depending on the rank of the bundle and the genus of the surface.


In this paper we extend the Daskalopoulos and Uhlenbeck approach to the moduli space of $\alpha$-stable representations of a quiver with fixed dimension vector, as introduced by King in \cite{King94} (see Section \ref{sec:background} for more information). The difference between the approach of this paper and Kirwan's approach of \cite{Kirwan86-2} is twofold: (a) Kirwan's results do not apply \emph{a priori} to the non-compact spaces $\Rep(Q, {\bf v})$ under consideration here, one first needs to prove a gradient flow convergence result for the norm-square of the moment map on $\Rep(Q, {\bf v})$ (which was done in \cite{haradawilkin}) as well as on the blow-ups of $\Rep(Q, {\bf v})$ required to desingularise the quotient (the proof of which has not yet been carried out), and (b) in this paper the dimensional bound $d_{min}$ can be computed directly by studying the different possible destabilising subrepresentations, whereas in Kirwan's method one needs to also study the desingularised spaces, the calculations for which become much more complicated as the dimension vector ${\bf v}$ increases in size.  Of course, it should be mentioned that Kirwan's methods do apply to a general class of quotients of compact K\"ahler manifolds and appear as a corollary of her results on the intersection cohomology of K\"ahler quotients, while the scope of this paper is limited to computing homotopy groups of K\"ahler quotients of $\Rep(Q, {\bf v})$.

The main result of this paper is a computation of the low-dimensional homotopy groups of the spaces $\Rep(Q, {\bf v})^{\alpha-st}$ and $\mathcal{M}_\alpha(Q, {\bf v})^{st}$, up to a dimension that depends on the underlying quiver, the choice of dimension vector, and the stability parameter $\alpha$. The following is the main theorem.

\begin{theorem}\label{thm:intro-homotopy-trivial}[Theorem \ref{thm:homotopy-trivial}]
Let $Q$ be a quiver, $\Rep(Q, {\bf v})$ be the space of representations with fixed dimension vector ${\bf v}$, let $\alpha \in \mathfrak{g}_{\bf v}^*$ be a stability parameter, and let $d_{min}(Q, {\bf v}, \alpha)$ be the associated minimal dimension as defined in Definition \ref{def:minimal-dimension}. Then $\pi_n( \Rep(Q, {\bf v})^{\alpha-st}) = 0$ for all $n < d_{min}(Q, {\bf v}, \alpha) - 1$.
\end{theorem}

\begin{corollary}\label{cor:intro-homotopy-moduli-space}[Corollary \ref{cor:homotopy-moduli-space}]
Let $PG_{\bf v} = G_{\bf v} /U(1)$ be the quotient of $G_{\bf v}$ by the subgroup consisting of scalar multiples of the identity. Then $\pi_n \left(\mathcal{M}_\alpha(Q, {\bf v})^{st} \right) \cong \pi_{n-1} (PG_{\bf v})$ for all $n < d_{min} (Q, {\bf v}, \alpha) - 1$.
\end{corollary}


As an application of the theory, some examples of interest are explicitly worked out in Section \ref{sec:examples}. The first example is a space that contains the moduli space of framed instantons on $S^4$, via the ADHM construction of \cite{ADHM78}, and while the ADHM moduli space itself (which is a hyperk\"ahler quotient, or alternately a K\"ahler quotient of a singular space) is beyond the scope of this paper, one can view these results as a first step towards studying the topology of this space. Again, it should be mentioned that in \cite{Kirwan92}, Kirwan takes a different approach to the problem of studying the ADHM moduli space, by first desingularising the hyperk\"ahler quotient of the space of semistable points, and then using general results (first derived in \cite{Kirwan86-3}) to describe the effect on low-dimensional homotopy groups of removing the strata corresponding to the singular points in the original moduli space.

The second example is the space of \emph{non-degenerate} polygons in $\R^3$ (i.e. those polygons that do not lie in a line) modulo rigid motions of $\R^3$. In Section \ref{subsec:polygon} we compute the homotopy groups of the moduli space of non-degenerate polygons up to dimension $2$ under certain conditions on the stability parameter.


Ultimately one would like to apply these methods to hyperk\"ahler quotients, in particular the Nakajima quiver varieties from \cite{Nakajima94}, and the moduli space of stable Higgs bundles of rank $n$ and degree $d$ over a compact Riemann surface from \cite{Hitchin87}. These can be viewed in two different ways: As a hyperk\"ahler quotient $\mu_1^{-1}(0) \cap \mu_2^{-1}(0) \cap \mu_3^{-1}(0) / G$, or as a GIT quotient of a singular space $\left( \mu_2^{-1}(0) \cap \mu_3^{-1}(0) \right) / \! / G^\C$. In the first case there is no stability condition that implies existence of common zeros of the three moment maps $\mu_1, \mu_2, \mu_3$, and in the second case the space $\mu_2^{-1}(0) \cap \mu_3^{-1}(0)$ is singular, and so \emph{a priori} the methods of this paper do not apply to this situation. It would be a very interesting problem to extend these techniques to stratified singular spaces in order to study hyperk\"ahler quotients such as those described above.

This paper is organised as follows: Section \ref{sec:background} contains the background material and notational conventions used in this paper, Section \ref{sec:transversality} contains the proof of the main theorem, and the explicit examples are worked out in Section \ref{sec:examples}.

\section{Background and notation}\label{sec:background}


A \emph{quiver} $Q$ is a directed graph, with vertices $\mathcal{I}$, edges $\mathcal{E}$ and head/tail maps $h, t : \mathcal{I} \rightarrow \mathcal{E}$ that define the orientation of the edges. A \emph{complex representation} of a quiver consists of a choice of complex vector spaces $\{ V_j \}_{j \in \mathcal{I}}$ with a fixed Hermitian structure, and homomorphisms $\{ A_a : V_{t(a)} \rightarrow V_{h(a)}\}_{a \in \mathcal{E}}$. The \emph{dimension vector} of a representation is the vector ${\bf v} = \left( \dim V_j \right)_{j \in \mathcal{I}}$. The space of representations with fixed dimension vector is denoted
\begin{equation}\label{eqn:def-representations}
\Rep(Q, {\bf v}) = \bigoplus_{a \in \mathcal{E}} \Hom(V_{t(a)}, V_{h(a)}),
\end{equation}
and the vector space consisting of the direct sum of all the vector spaces at each vertex is denoted
\begin{equation*}
\Vect(Q, {\bf v}) = \bigoplus_{j \in \mathcal{I}} V_j .
\end{equation*}
 
A \emph{subrepresentation} of a representation $A \in \Rep(Q, {\bf v})$ consists of vector spaces $\{ V_j' \subseteq V_j \}_{j \in \mathcal{I}}$ such that $A_a (V_{t(a)}') \subseteq V_{h(a)}'$ for all $a \in \mathcal{E}$, and homomorphisms $\{ A_a' : V_{t(a)}' \rightarrow V_{h(a)}' \}_{a \in \mathcal{E}}$ such that each $A_a'$ is the restriction of $A_a$ to $V_{t(a)}'$. 

Given such a subrepresentation, for each vertex $j \in \mathcal{I}$ there is an associated projection $\pi_j : V_j \rightarrow V_j'$ that is defined using orthogonal projection with respect to the Hermitian structure on $V_j$. If the dimension vector of the subrepresentation is ${\bf v'} = ( \dim V_j' )_{j \in \mathcal{I}}$ then we write ${\bf v'} \leq {\bf v}$, the projection $\pi : \Vect(Q, {\bf v}) \rightarrow \Vect(Q, {\bf v'})$ is used to denote the vector spaces in the subrepresentation, and $A' \in \Rep(Q, {\bf v'})$ is used to denote the induced homomorphisms $A_a' = \pi \circ A_a \circ \pi \in \Hom(V_{t(a)}', V_{h(a)}')$. 


Given a fixed Hermitian structure on each vector space $V_j$, the group $G_{\bf v} = \times_{j \in \mathcal{I}} {\sf U}(V_j)$ acts on $\Rep(Q, {\bf v})$ by
\begin{equation}\label{eqn:def-group-action}
g \cdot \{ A_a \}_{a \in \mathcal{E}} = \{ g_{h(a)} A_a g_{t(a)}^{-1} \}_{a \in \mathcal{E}} ,
\end{equation}
where $g = \left( g_j \right)_{j \in \mathcal{I}} \in \times_{j \in \mathcal{I}} {\sf U}(V_j)$. Let $\mathfrak{g}_{\bf v}$ denote the Lie algebra of $G_{\bf v}$. The Hermitian structure on each vector space $V_j$ gives an identification 
\begin{equation*}
\mathfrak{g}_{\bf v} = \bigoplus_{j \in \mathcal{I}} \mathfrak{u}(V_j) \cong \bigoplus_{j \in \mathcal{I}} \mathfrak{u}(V_j)^* = \mathfrak{g}_{\bf v}^*,
\end{equation*} 
and with respect to this Hermitian structure we define the adjoint $A_a^*$ of each homomorphism $A_a : V_{t(a)} \rightarrow V_{h(a)}$.

The action of $G_{\bf v}$ preserves the symplectic structure (in fact the K\"ahler structure) on $\Rep(Q, {\bf v})$ induced by the Hermitian structure on each $V_j$, and a moment map for this action is a map $\mu : \Rep(Q, {\bf v}) \rightarrow \mathfrak{g}_{\bf v}^* = \bigoplus_{j \in \mathcal{I}} \mathfrak{u}(V_j)^*$, where the component corresponding to $\mathfrak{u}(V_j)^*$ is
\begin{equation}\label{eqn:moment-map-component-def}
\mu_j(A) = \sqrt{-1} \sum_{\substack{a \in \mathcal{E} \\ h(a) = j}} A_a A_a^* - \sqrt{-1} \sum_{\substack{a \in \mathcal{E} \\ t(a) = j}} A_a^* A_a .
\end{equation}
For convenience, in the rest of the paper the notation
\begin{equation}\label{eqn:moment-map-def}
\mu(A) = \sqrt{-1} \sum_{a \in \mathcal{E}} [A_a, A_a^*] \in \mathfrak{g}_{\bf v}^* \subset \End \left( \Vect(Q, {\bf v}) \right)
\end{equation}
is used, where each $A_a \in \Hom(V_{t(a)}, V_{h(a)}) \subseteq \End \left( \Vect(Q, {\bf v}) \right)$, and the commutator of two elements $X, Y \in \End \left( \Vect(Q, {\bf v}) \right)$ is denoted $[X, Y]$.

A central element $\alpha$ of $\mathfrak{g}_{\bf v}^*$ consists of a choice of a real number $\alpha_j$ at each vertex $j \in \mathcal{I}$, denoted 
\begin{equation*}
\alpha = \left( \sqrt{-1} \alpha_j \id_{V_j} \right)_{j \in \mathcal{I}} .
\end{equation*}
Corollary 6.2 of \cite{King94} identifies the K\"ahler quotient 
\begin{equation}\label{eqn:Kahler-quotient}
\mathcal{M}_\alpha(Q, {\bf v}) : = \mu^{-1}(\alpha) / G_{\bf v} 
\end{equation}
with the \emph{moduli space of $\alpha$-semistable representations}, $\Rep(Q, {\bf v}) / \negthinspace /_\alpha G_{\bf v}^\C$, defined using GIT (see below for more details). This space depends on the quiver $Q$, the dimension vector ${\bf v}$, and the choice of $\alpha \in \mathfrak{g}_{\bf v}^*$, and note that since $\mu(A)$ consists of commutators then the quotient is empty if $\sum_{j \in \mathcal{I}} \alpha_j \dim_\C V_j \neq 0$.


The construction of King in \cite{King94} shows that when $\alpha$ is integral then this quotient can be described via Geometric Invariant Theory as follows.  Let $\alpha = \left( \sqrt{-1} \alpha_j \id_{V_j} \right)_{j \in \mathcal{I}}$, with each $\alpha_j \in \ZBbb$, and let $\displaystyle{G_{\bf v}^\C = \times_{j \in \mathcal{I}} \GL(V_j, \C)}$ denote the complexification of $G$. The action of the group $G_{\bf v}^\C$ lifts to an action on the trivial line bundle over $\Rep(Q, {\bf v})$ by
\begin{equation*}
g \cdot \left( \left\{ A_a \right\}_{a \in \mathcal{E}}, \xi \right) = \left( \left\{ g_{h(a)} A_a g_{t(a)}^{-1} \right\}, \chi(g) \xi \right), \, \text{for} \, \xi \in \C,
\end{equation*}
where the character $\chi$ is given by
\begin{equation*}
\chi(g) = \prod_{j \in \mathcal{I}} \det(g_j)^{-\alpha_j} .
\end{equation*}
A representation $A \in \Rep(Q, {\bf v})$ is \emph{$\alpha$-semistable} if (given some nonzero $\xi \in \C$) the closure of the $G_{\bf v}^\C$-orbit of $(A, \xi)$ in $\Rep(Q, {\bf v}) \times \C$ does not intersect the zero section, i.e.
\begin{equation*}
\overline{ G_{\bf v}^\C \cdot (A, \xi) } \cap \left( \Rep(Q, {\bf v}) \times \{ 0 \} \right) = \emptyset .
\end{equation*}
A representation $A \in \Rep(Q, {\bf v})$ is \emph{$\alpha$-stable} if (given some nonzero $\xi \in \C$) the $G_{\bf v}^\C$-orbit of $(A, \xi)$ in $\Rep(Q, {\bf v}) \times \C$ is closed.

Proposition 3.1 in \cite{King94} shows that the definitions given above are equivalent to the following slope-stability conditions (analogous to those for holomorphic bundles over a K\"ahler manifold). Given a subrepresentation with dimension vector ${\bf v'} = \left( v_j \right)_{j \in \mathcal{I}}$, the \emph{rank} and \emph{$\alpha$-degree} of the subrepresentation are given by
\begin{align}
\rank(Q, {\bf v'}) & = \sum_{j \in \mathcal{I}} v_j \label{eqn:def-rank} \\
\deg_\alpha(Q, {\bf v'}) & = \sum_{j \in \mathcal{I}} - \alpha_j v_j  \label{eqn:def-degree}
\end{align}
The \emph{$\alpha$-slope} is then defined to be the quotient $\slope_\alpha(Q, {\bf v'}) = \deg_\alpha(Q, {\bf v'}) / \rank(Q, {\bf v'})$.
\begin{lemma}[Proposition 3.1 of \cite{King94}]\label{lem:slope-stability}
Let $Q$ be a quiver, ${\bf v}$ a dimension vector and $\alpha$ a stability parameter satisfying $\deg_\alpha(Q, {\bf v}) = 0$. A representation $A \in \Rep(Q, {\bf v})$ is $\alpha$-stable (resp. $\alpha$-semistable) if and only if every subrepresentation satisfies
\begin{equation}\label{eqn:def-slope-stability}
\slope_\alpha(Q, {\bf v'}) < \slope_\alpha(Q, {\bf v}) \quad \left(\text{resp.} \, \slope_\alpha(Q, {\bf v'}) \leq \slope_\alpha(Q, {\bf v}) \, \right) .
\end{equation} 
\end{lemma}

The space of $\alpha$-stable (resp. $\alpha$-semistable) representations is denoted $\Rep(Q, {\bf v})^{\alpha-st}$ (resp. $\Rep(Q, {\bf v})^{\alpha-ss}$), and the GIT quotient is defined to be
\begin{equation}\label{eqn:GIT-quotient}
\mathcal{M}_\alpha(Q, {\bf v}) = \Rep(Q, {\bf v}) / \negthickspace /_\alpha G_{\bf v}^\C :=  \Rep(Q, {\bf v})^{\alpha-ss} / \negthickspace / G_{\bf v}^\C,
\end{equation}
where the quotient $/ \negthickspace /$ identifies orbits whose closures intersect.


It follows from Theorem 6.1 of \cite{King94} that the K\"ahler quotient defined in \eqref{eqn:Kahler-quotient} can be identified with the GIT quotient defined in \eqref{eqn:GIT-quotient} (see also \cite{haradawilkin} for an alternative proof using the gradient flow of $\| \mu - \alpha \|^2$).



Reineke in \cite{Reineke03} computes Betti numbers of $\mathcal{M}_\alpha(Q, {\bf v})^{st}$ for choices of parameter and dimension vector for which the sets $\Rep(Q, {\bf v})^{\alpha-st}$ and $\Rep(Q, {\bf v})^{\alpha-ss}$ co-incide. Following on from this, the results of \cite{haradawilkin} show that there is a surjective map (known as the \emph{Kirwan map}) $\kappa : H_{G_{\bf v}}^*(\Rep(Q, {\bf v})) \twoheadrightarrow H_{G_{\bf v}}^*( \Rep(Q, {\bf v})^{\alpha-ss})$ and that the Betti numbers can also be computed using the equivariant Morse theory of Atiyah and Bott \cite{AtiyahBott83} and Kirwan \cite{Kirwan84}. The Morse theory only gives information about the topology of $\Rep(Q, {\bf v})^{\alpha-ss}$ and the purpose of this paper is to understand the low-dimensional topology of $\mathcal{M}_\alpha(Q, {\bf v})^{\alpha-st}$ regardless of whether $\Rep(Q, {\bf v})^{\alpha-st}$ and $\Rep(Q, {\bf v})^{\alpha-ss}$ co-incide.

If the sets $\Rep(Q, {\bf v})^{\alpha-st}$ and $\Rep(Q, {\bf v})^{\alpha-ss}$ do not co-incide then the quotient may have serious singularities (see for example \cite{Kirwan85}). In these cases there is a smooth structure on the quotient of the subset of stable representations
\begin{equation}
\mathcal{M}_\alpha(Q, {\bf v})^{st} = \Rep(Q, {\bf v})^{\alpha-st} / G_{\bf v}^\C .
\end{equation}

The situation in this paper is an extension to the case of quivers of the infinite-dimensional methods for studying the moduli space of stable holomorphic bundles over a compact Riemann surface used in \cite{Daskal92} and \cite{DaskalUhlenbeck95}. For  rank $2$ bundles Daskalopoulos used an algebraic stratification of the space of semistable bundles to obtain results about the low-dimensional topology of the space of stable bundles (see Section 7 of \cite{Daskal92} for more details). As noted on pp744-745 of \cite{Daskal92}, for higher rank bundles a stratification argument quickly becomes intractable, and in \cite{DaskalUhlenbeck95} a more general approach using transversality is adopted. In this paper we follow the technique of \cite{DaskalUhlenbeck95}, and use a transversality argument to obtain the results of Theorem \ref{thm:homotopy-trivial}.


Given a representation $A \in \Rep(Q, {\bf v})$, the \emph{infinitesimal action} of $G_{\bf v}^\C$ at $A$ is given by
\begin{align}\label{eqn:inf-action-def}
\begin{split}
\rho_A : \mathfrak{g}_{\bf v}^\C & \rightarrow T_A \Rep(Q, {\bf v}) \\
 u & \mapsto \left. \frac{d}{dt} \right|_{t=0} e^{tu} \cdot A = \left( [u, A_a] \right)_{a \in \mathcal{E}} .
\end{split}
\end{align}

Given $A \in \Rep(Q, {\bf v})^{\alpha-st}$, the infinitesimal deformations of the moduli space $\mathcal{M}_\alpha(Q, {\bf v})^{st}$ at the equivalence class containing $A$ can be viewed in terms of the cohomology of a complex
\begin{equation}\label{eqn:deformation-complex-def}
\mathfrak{g}_{\bf v}^\C \stackrel{\rho_A}{\rightarrow} T_A \Rep(Q, {\bf v}) .
\end{equation}

\begin{remark}\label{rem:stable-cohomological-interpretation}
This complex can be interpreted in terms of the cohomology of the quiver representation $V_A$ associated to $A$, via the following exact sequence, which first appeared in Section 2.1 of \cite{Ringel76} (I thank the referee for this reference)
\begin{multline}
0 \rightarrow \Hom(V_A, V_A) \cong \ker \rho_A \rightarrow \mathfrak{g}_\C \stackrel{\rho_A}{\rightarrow}  T_A \Rep(Q, {\bf v})
 \rightarrow \coker \rho_A \cong \Ext^1(V_A, V_A) \rightarrow 0 .
\end{multline}
\end{remark}

The tangent space to the moduli space is isomorphic to $T_A \Rep(Q, {\bf v}) / \im \rho_A = \coker \rho_A\cong \Ext^1(V_A, V_A)$. Since an $\alpha$-stable point has stabiliser consisting of scalar multiples of the identity in $G_{\bf v}^\C$, then $\dim_\C \ker \rho_A = 1$, and so
\begin{align*}
\dim_\C \coker \rho_A - 1& = \dim_\C \coker \rho_A - \dim_\C \ker \rho_A \\
 & = \dim_\C T_A \Rep(Q, {\bf v}) - \dim_\C \mathfrak{g}_{\bf v}^\C \\
\Leftrightarrow \quad \dim_\C \coker \rho_A & = \dim_\C T_A \Rep(Q, {\bf v}) - \dim_\C \mathfrak{g}_{\bf v}^\C + 1.
\end{align*}

When $A \in \Rep(Q, {\bf v})$ is $\alpha$-unstable, then one can still define the complex \eqref{eqn:deformation-complex-def}, however $\ker \rho_A$ may now be non-trivial. A particular case of interest is when $A$ preserves the vector spaces $\Vect(Q, {\bf v'})$ of a given subrepresentation with associated orthogonal projection $\pi$. In this case one can view $\pi$ as an endomorphism $\pi \in \mathfrak{g}_{\bf v}^\C = \bigoplus_{j \in \mathcal{I}} \End(V_j)$, and $\rho_A (\pi) = \left( [A_a, \pi] \right)_{a \in \mathcal{E}}$. If $\pi : \Vect(Q, {\bf v}) \rightarrow \Vect(Q, {\bf v'})$ is a subrepresentation of $A$ (i.e. $A_a \left(\pi (V_{t(a)}) \right) \subseteq \pi(V_{h(a)})$ for each $a \in \mathcal{E}$), then $\pi \circ A_a \circ \pi = A_a \circ \pi = A_a \circ \pi \circ \pi$, or alternately $[A_a, \pi] \circ \pi = 0$ for each $a \in \mathcal{E}$. Conversely, if $[A_a, \pi] \circ \pi = 0$ for each $a \in \mathcal{E}$, then $A_a$ preserves the vector spaces in the image of the projection $\pi$. Therefore we have proved the following lemma.

\begin{lemma}
$\pi : \Vect(Q, {\bf v}) \rightarrow \Vect(Q, {\bf v'})$ is a subrepresentation of $A \in \Rep(Q, {\bf v})$ iff $\rho_A (\pi) \circ \pi = 0$.
\end{lemma}

Given a subrepresentation $\pi : \Vect(Q, {\bf v}) \rightarrow \Vect(Q, {\bf v'})$ and the orthogonal complement $\pi^\perp := (\id - \pi) : \Vect(Q, {\bf v}) \rightarrow \Vect(Q, {\bf v}')^\perp$ we define the following subspaces of $\mathfrak{g}_{\bf v}^\C$ and $T_A \Rep(Q, {\bf v})$.
\begin{align*}
\Hom^0 (\pi, \pi^\perp) & : = \bigoplus_{j \in \mathcal{I}} \Hom \left( \pi(V_j), \pi^\perp(V_j) \right) \subseteq \mathfrak{g}_{\bf v}^\C \\
\Hom^1 (\pi, \pi^\perp) & : = \bigoplus_{a \in \mathcal{E}} \Hom \left( \pi(V_{t(a)}), \pi^\perp(V_{h(a)}) \right) \subseteq \Rep(Q, {\bf v})
\end{align*}
Moreover, if $A$ preserves $\pi$ then $\rho_A : \mathfrak{g}_\C \rightarrow T_A \Rep(Q, {\bf v})$ induces a map
\begin{equation} \label{eqn:rho_pi}
\Hom^0 (\pi, \pi^\perp) \stackrel{\rho_\pi}{\rightarrow} \Hom^1 (\pi, \pi^\perp),
\end{equation}
where we identify $\Hom^1(\pi, \pi^\perp) \subseteq \Rep(Q, {\bf v})$ with a quotient of the tangent space $T_A \Rep(Q, {\bf v})$. If $A$ preserves both $\pi$ and $\pi^\perp$, then we say that the representation \emph{splits}. This viewpoint will become useful in the next section, where we use this to define a simple condition for a representation to be unstable.

\begin{remark}\label{rem:unstable-cohomological-interpretation}
The complex \eqref{eqn:rho_pi} also has a cohomological interpretation similar to that in Remark \ref{rem:stable-cohomological-interpretation} via the following complex
\begin{equation}
0 \rightarrow \Hom(V_\pi, V / V_\pi) \rightarrow \Hom^0(\pi, \pi^\perp)\stackrel{\rho_\pi}{\rightarrow}   \Hom^1(\pi, \pi^\perp) \rightarrow  \Ext(V_\pi, V / V_\pi) \rightarrow 0 ,
\end{equation}
where $\Hom(V_\pi, V / V_\pi) \cong \ker \rho_\pi$ and $\Ext(V_\pi, V / V_\pi) \cong \coker \rho_\pi$. The Euler characteristic of the subrepresentation (defined below) is then fixed independently of the choice of representations $V_\pi$ and $V / V_\pi$
\begin{align}\label{eqn:Euler-characteristic}
\begin{split}
\chi(V_\pi, V / V_\pi) & := \dim_\C \Hom (V_\pi, V / V_\pi) - \dim_\C \Ext^1 (V_\pi, V / V_\pi) \\
 & = \dim_\C \Hom^0(\pi, \pi^\perp) - \dim_\C \Hom^1(\pi, \pi^\perp) .
\end{split}
\end{align}
When the representation splits, then $V / V_\pi \cong V_{\pi^\perp}$, and so we denote
\begin{equation*}
\chi(\pi, \pi^\perp) : = \chi(V_\pi, V_{\pi^\perp}) .
\end{equation*}
\end{remark}

\section{Transversality results}\label{sec:transversality}

The purpose of this section is to prove Theorem \ref{thm:homotopy-trivial}, which shows that the homotopy groups of the space $\Rep(Q, {\bf v})^{\alpha-st}$ are trivial up to a dimensional bound that depends on the quiver $Q$, the dimension vector ${\bf v}$ and the parameter $\alpha$. The key idea behind this is to show that for small values of $n$ there exists a based homotopy equivalence $F$ between a given map $s : S^n \rightarrow \Rep(Q, {\bf v})^{\alpha-st}$ and the constant map, such that the image of $F$ is contained in $\Rep(Q, {\bf v})^{\alpha-st}$. Equivalently, we show that the image of $F$ does not contain any representations that preserve a sub-representation violating the slope-stability condition of Lemma \ref{lem:slope-stability}.

Given a dimension vector ${\bf v'} \leq {\bf v}$, the set of all possible subrepresentations $\pi : \Vect(Q, {\bf v}) \rightarrow \Vect(Q, {\bf v'})$ is given by
\begin{equation}\label{eqn:def-grassmannian}
\Gr({\bf v'}, {\bf v}) : = \prod_{j \in \mathcal{I}} \Gr(v_j', v_j) .
\end{equation}
For $\pi \in \Gr({\bf v'}, {\bf v})$, the value of $\dim_\R \Hom^1 (\pi, \pi^\perp) - \dim_\R \Hom^0 (\pi, \pi^\perp)$ depends only on the dimension vector ${\bf v'}$. 

\begin{definition}\label{def:minimal-dimension}
The \emph{minimal dimension} of the triple $(Q, {\bf v}, \alpha)$ is 
\begin{align}\label{eqn:minimal-dimension}
\begin{split}
d_{min}(Q, {\bf v}, \alpha) & : = \min_{\substack{{\bf v'} < {\bf v} \\ \slope_\alpha(Q, {\bf v'}) \geq \slope_\alpha(Q, {\bf v}) }} \dim_\R  \Hom^1 (\pi, \pi^\perp) - \dim_\R \Hom^0 (\pi, \pi^\perp) \\
 & = \min_{\substack{{\bf v'} < {\bf v} \\ \slope_\alpha(Q, {\bf v'}) \geq \slope_\alpha(Q, {\bf v}) }} (-2) \chi(\pi, \pi^\perp).
\end{split}
\end{align}
\end{definition}

\begin{remark}
Theorem \ref{thm:homotopy-trivial} and Corollary \ref{cor:homotopy-moduli-space} show that the value of $d_{min}(Q, {\bf v}, \alpha) - 2$ gives an upper bound on the dimension of the homotopy groups of $\mathcal{M}(Q, {\bf v})^{\alpha-st}$ that can be computed using the methods of this paper. In Section \ref{sec:examples} we compute $d_{min}(Q, {\bf v}, \alpha)$ explicitly for each example studied, and it would be an extremely interesting problem to provide some general combinatorial bounds for $d_{min}(Q, {\bf v}, \alpha)$ in terms of the quiver and dimension vector. Results of Schofield in \cite{Schofield92} provide a starting point in this direction, as they give some conditions for which (in the notation of this paper) $\Ext(V_\pi, V_{\pi^\perp})$ vanishes generally, which relates to $d_{min}$ via \eqref{eqn:Euler-characteristic}. 
\end{remark}

Let $T \rightarrow \Gr({\bf v'}, {\bf v})$ denote the tautological bundle with fibre over $\pi \in \Gr({\bf v'}, {\bf v})$ given by $\Vect(Q, {\bf v'})$, a sub-bundle of the trivial $\Vect(Q, {\bf v})$ bundle over $\Gr({\bf v'}, {\bf v})$. Similarly, let $W^0 \rightarrow \Gr({\bf v'}, {\bf v})$ be the tautological sub-bundle of the trivial bundle $\Gr({\bf v'}, {\bf v}) \times \mathfrak{g}_{\bf v}^\C$, where the fibres of $W^0$ are the subspaces $\Hom^0(\pi, \pi^\perp) \subseteq \mathfrak{g}_{\bf v}^\C$, and let $W^1 \rightarrow \Gr({ \bf v'}, {\bf v})$ be the tautological sub-bundle of the trivial bundle $\Gr({\bf v'}, {\bf v}) \times \Rep(Q, {\bf v})$, where the fibres of $W^1$ are the subspaces $\Hom^1(\pi, \pi^\perp) \subseteq \Rep(Q, {\bf v})$.

Given a compact manifold $X$, and basepoints $x_0 \in X$ and $A_0 \in \Rep(Q, {\bf v})$, let $C^q \left(X, \Rep(Q, {\bf v}) \right)$ be the set of $C^q$ maps from $X$ to $\Rep(Q, {\bf v})$, this has the structure of a $C^q$ manifold. Let $\mathcal{C}(X) \subset C^q \left(X, \Rep(Q, {\bf v}) \right)$ be the subset of basepoint-preserving $C^q$ maps from $X$ to $\Rep(Q, {\bf v})$. 

\begin{lemma}
$\mathcal{C}(X)$ is a submanifold of $C^q \left( X, \Rep(Q, {\bf v}) \right)$, and the tangent space at a point $s \in \mathcal{C}(X)$ is
\begin{equation}
T_s \mathcal{C}(X) = \left\{ \delta s : X \rightarrow \Rep(Q, {\bf v}) \, : \, \delta s \in C^q, \delta s (x_0) = 0 \right\} .
\end{equation}
\end{lemma}

\begin{proof}
Define $L : C^q \left(X, \Rep(Q, {\bf v}) \right) \rightarrow \Rep(Q, {\bf v})$ by $L(s) = s(x_0)$. Then $\mathcal{C}(X) = L^{-1} (A_0)$, and so by the implicit function theorem, the proof reduces to showing that $A_0$ is a regular value of $L$. The derivative of $L$ at $s \in \mathcal{C}(X)$ is $\delta L_s (\delta s) = \delta s (x_0)$, and the Whitney extension theorem (cf. Appendix A in \cite{AbrahamRobbin67}) shows that given a prescribed derivative $\delta s$ at $x_0$ there exists $s \in C^q \left(X, \Rep(Q, {\bf v}) \right)$ with $s(x_0) = A_0$ and the derivative of $s$ at $x_0$ equal to $\delta s$. Therefore $\delta L$ is surjective and $A_0$ is a regular value of $L$.
\end{proof}

The space of basepoint-preserving homotopy equivalences between maps in $\mathcal{C}(X)$  is defined as follows. Let $\mathcal{S}(X) \subset C^q \left(X \times [0,1],\Rep(Q, {\bf v}) \right)$ be the subset of maps $F : X \times [0,1] \rightarrow \Rep(Q, {\bf v})$ such that $F(x_0, t) = A_0$ for all $t \in [0,1]$. 
\begin{lemma}
The space $\mathcal{S}(X)$ has the structure of a $C^q$ manifold, and the tangent space at $F \in \mathcal{S}(X)$ is
\begin{equation}
T_F \mathcal{S}(X) = \left\{ \delta F : X \times [0,1] \rightarrow \Rep(Q, {\bf v}) \, : \, \delta F \in C^q, \delta F(x_0, t) = 0 \, \forall t \in [0,1] \right\} .
\end{equation}
\end{lemma}

The proof of this statement is similar to that of the previous lemma (which uses the implicit function theorem and Whitney extension theorem), and so it is omitted.

Given $s \in \mathcal{C}(X)$, the space of basepoint-preserving homotopy equivalences between $s$ and the constant map is
\begin{equation}
\mathcal{S}(s) = \left\{ F \in \mathcal{S}(X) \, : \, F(x, 0) = s(x), F(x, 1) = A_0, F(x_0, t) = A_0 \right\} .
\end{equation}
Since $\Rep(Q, {\bf v})$ is a vector space then $\mathcal{S}(s)$ is clearly non-empty.

\begin{lemma}
$\mathcal{S}(s)$ is a $C^q$ submanifold of $\mathcal{S}(X)$, and the tangent space of $\mathcal{S}(s)$ at $F$ is 
\begin{multline}\label{eqn:tangent-based-homotopy}
T_F \mathcal{S}(s) = \{ \delta F : X \times [0,1] \rightarrow \Rep(Q, {\bf v}) \, : \, \\ \delta F \in C^q, \delta F(x,0) = \delta F(x, 1) = \delta F(x_0, t) = 0, \forall t \in [0,1], x \in X \} .
\end{multline}
\end{lemma}

The proof of this lemma is similar to that of the previous lemma (see also the proof of Lemma 2.1 in \cite{DaskalUhlenbeck95}), and so it is omitted.

Now fix $X = S^n$ and $s : X \rightarrow \Rep(Q, {\bf v})^{\alpha-st}$. For notation let $\mathcal{C} = \mathcal{C}(X)$, $\mathcal{S} = \mathcal{S}(s)$ and $Y = S^n \times [0,1]$. Given a dimension vector ${\bf v'} \leq {\bf v}$, define the map 
\begin{align*}
\mathcal{D} : \Gr({\bf v'}, {\bf v}) \times \mathcal{S} \times Y  & \rightarrow W^1 \subseteq \Gr({\bf v'}, {\bf v}) \times \Rep(Q, {\bf v}) \\
\mathcal{D} (\pi, F, y) & = \left( \pi, \pi^\perp \circ \rho_{F(y)}(\pi) \circ \pi \right) 
\end{align*}
Since the component of $\rho_{F(y)}(\pi)$ along each edge $a \in \mathcal{E}$ is given by $\left[ F(y)_a, \pi \right]$, then for simplicity of notation we write $\rho_{F(y)}(\pi) = \left[ F(y), \pi \right]$, where the meaning is that the component of $\left[ F(y), \pi \right]$ corresponding to the edge $a \in \mathcal{E}$ is $\left[ F(y)_a, \pi \right]$. For a given map $F \in \mathcal{S}(s)$, define $\mathcal{D}^F (\pi, y) = \mathcal{D}(\pi, F, y)$.

\begin{remark}\label{rem:simplify-D}

\begin{enumerate}
\item Since $\pi^\perp \circ \pi = 0$ and $\pi \circ \pi = \pi$, then we can write 
\begin{equation*}
\mathcal{D} (\pi, F, y) = \left( \pi, \pi^\perp \circ [F(y), \pi] \circ \pi \right) = \left( \pi, \pi^\perp \circ F(y) \circ \pi \right) .
\end{equation*}

\item Let $\mathcal{O}_{W^1}$ denote the zero section of $W^1$. Then $F(y) \in \Rep(Q, {\bf v})$ preserves the subspaces defined by $\pi$ iff $\pi^\perp \circ F(y) \circ \pi = 0$, which by the previous statement is equivalent to saying that $F(y)$ preserves $\pi$ iff $(\pi, F, y) \in \mathcal{D}^{-1} (\mathcal{O}_{W^1})$. 

\item The same computation as above shows that 
\begin{equation*}
F(y) \, \, \text{preserves} \, \,  \pi \quad \Leftrightarrow \quad \pi^\perp \circ [F(y), \pi] = 0 \quad \Leftrightarrow \quad [F(y), \pi] \circ \pi = 0 .
\end{equation*}
\end{enumerate}
\end{remark}

\begin{proposition}\label{prop:transversality-conditions}

\begin{enumerate}

\item \label{item:fibre} $\mathcal{D} : \Gr({\bf v'}, {\bf v}) \times \mathcal{S} \times Y \rightarrow W^1$ is a $C^q$ map of fibre bundles over $\Gr({\bf v'}, {\bf v})$. 

\item \label{item:transverse} $\mathcal{D}$ is transverse to the zero section $\mathcal{O}_{W^1}$ of $W^1$.

\item \label{item:index} Let $W_\pi^1$ be the fibre of $W^1$ over $\pi \in \Gr({\bf v'}, {\bf v})$. Then 
\begin{align}
\codim_{W^1} \mathcal{O}_{W^1} & = \dim_\R W_\pi^1 = \dim_\R \Hom^1 (\pi, \pi^\perp) \label{eqn:codimension-condition} \\
\dim_\R \Gr({\bf v'}, {\bf v}) \times Y & = \dim_\R \Hom^0 (\pi, \pi^\perp) + n + 1 \label{eqn:grass-dim}
\end{align}
\end{enumerate}
\end{proposition}

\begin{proof}

It is clear that $\mathcal{D}$ maps fibres of $\Gr({\bf v'}, {\bf v}) \times \mathcal{S} \times Y$ to fibres of $W^1$, and so the first statement follows from the fact that $F$ is $C^q$.

To see the second statement, a calculation shows that the derivative of $\mathcal{D}$ with respect to $F$ at a point $(\pi, F, y) \in \mathcal{D}^{-1} (\mathcal{O}_{W^1})$is given by
\begin{equation}
d_2 \mathcal{D}_{(\pi, F, y)} ( \delta F) = \pi^\perp \circ \delta F \circ \pi .
\end{equation}
Since $s : X \rightarrow \Rep(Q, {\bf v})^{\alpha-st}$, $A_0 \in \Rep(Q, {\bf v})^{\alpha-st}$, and $F \in \mathcal{S}(s)$, then $(\pi, F, y) \in \mathcal{D}^{-1} (\mathcal{O}_{W^1})$ implies that $y \notin X \times \{ 0 \}$, $y \notin X \times \{ 1 \}$ and $y \notin \{ x_0 \} \times [0,1]$. Therefore, from \eqref{eqn:tangent-based-homotopy}, $\delta F(y)$ can take any value in $\Rep(Q, {\bf v})$, and so the above equation shows that $d_2 \mathcal{D}_{(\pi, F, y)}$ is onto $\Hom^1(\pi, \pi^\perp)$, the fibre of $W^1$ at $\pi$.

In the proof of the third statement, equation \eqref{eqn:codimension-condition} follows by observing that $W_\pi^1 \cong \Hom^1(\pi, \pi^\perp)$. Equation \eqref{eqn:grass-dim} then follows from the definition of $\Gr({\bf v'}, {\bf v})$ in \eqref{eqn:def-grassmannian}.
\end{proof}

In the following Corollary, recall that a \emph{residual subset} is the countable intersection of open dense subsets.


\begin{corollary}\label{cor:non-intersecting}
If $q \geq 1$ and $\dim_\R \Hom^0(\pi, \pi^\perp) + n + 1 < \dim_\R \Hom^1 (\pi, \pi^\perp)$, then for $F$ in a residual subset of $\mathcal{S}$, the image $\mathcal{D}^F \left( \Gr({\bf v'}, {\bf v} ) \times S^n \times [0,1] \right)$ does not intersect $\mathcal{O}_{W^1}$.
\end{corollary}

\begin{proof}


For $(\pi, F, y) \in \mathcal{D}^{-1}(\mathcal{O}_{W^1})$, let $p : T_{\mathcal{D}(\pi, F, y)} W^1 \rightarrow T_0 W_\pi^1$ be the natural projection map. Since the manifolds $\Gr({\bf v'}, {\bf v})$ and $Y$ are finite-dimensional, as are the fibres $W_\pi^1$ of the bundle $W^1$, then $p \circ (\delta \mathcal{D}^F)_{(\pi, y)} : T_\pi \Gr({\bf v'}, {\bf v}) \times T_y Y \rightarrow T_0 W_\pi^1$ is Fredholm of index 
\begin{equation*}
\dim_\R \Gr({\bf v'}, {\bf v}) + \dim_\R Y - \dim_\R W_\pi^1 = \dim_\R \Hom^0(\pi, \pi^\perp) + n + 1 - \dim \Hom^1 (\pi, \pi^\perp) .
\end{equation*}

Therefore a combination of Proposition \ref{prop:transversality-conditions} together with the Daskalopoulos-Uhlenbeck version of the transversal density theorem (Theorems A2 and A4 in \cite{DaskalUhlenbeck95}) shows that the set 
\begin{equation*}
\mathcal{S}_{\mathcal{O}_{W^1}} = \left\{ F \in \mathcal{S}\, : \, \mathcal{D}^F \, \text{ is transversal to } \, \mathcal{O}_{W^1} \right\}
\end{equation*}
is residual in $\mathcal{S}$.


Since the dimension of the fibres of $W^1$ (or alternately, the codimension of $\mathcal{O}_{W^1}$ in $W^1$) is given by $\dim \Hom^1 (\pi, \pi^\perp)$, and by assumption we have $\dim \Gr({\bf v'}, {\bf v}) + n + 1 < \dim \Hom^1 (\pi, \pi^\perp)$, then the map $\mathcal{D}^F : \Gr( {\bf v'}, {\bf v}) \times S^{n} \times [0,1] \rightarrow W^1$ is transverse to the zero section $\mathcal{O}_{W^1}$  iff the image $\mathcal{D}^F \left( \Gr( {\bf v'}, {\bf v}) \times S^{n} \times [0,1]  \right)$ does not intersect $\mathcal{O}_{W^1}$.
\end{proof}

\begin{remark}
This corollary shows that there is a residual set of $C^1$ basepoint-preserving homotopy equivalences $F$ between $s : S^n \rightarrow \Rep(Q, {\bf v})^{\alpha-st}$ and the constant map, such that no representation in the image of $F$ preserves a sub-representation in $\Gr({\bf v'}, {\bf v})$. 
\end{remark}

\begin{theorem}\label{thm:homotopy-trivial}
If $n + 1 < d_{min} (Q, {\bf v}, \alpha)$ then $\pi_n( \Rep(Q, {\bf v})^{\alpha-st}) = 0$.
\end{theorem}

\begin{proof}

Fix a continuous map $s : S^n \rightarrow \Rep(Q, {\bf v})^{\alpha-st}$. We aim to show that there is a basepoint-preserving homotopy equivalence in $\Rep(Q, {\bf v})^{\alpha-st}$ between $s$ and the constant map. Firstly note that since the set of stable points $\Rep(Q, {\bf v})^{\alpha-st}$ is open in the vector space $\Rep(Q, {\bf v})$, then we can perturb $s$ (while remaining in $\Rep(Q, {\bf v})^{\alpha-st}$) so that $s \in C^1$ (cf. Theorem 10.21 of \cite{Lee03}).

Corollary \ref{cor:non-intersecting} then shows that for each ${\bf v'} < {\bf v}$ there is a residual set $\mathcal{S}(s)_{\bf v'}^0$ of homotopy equivalences $F$ from $s$ to the constant map such that no representation in the image of $F$ preserves a subrepresentation in $\Gr({\bf v'}, {\bf v})$.

There are a finite number of choices of ${\bf v'} < {\bf v}$ such that $\slope_\alpha(Q, {\bf v'}) \geq \slope_\alpha(Q, {\bf v})$, and so the set
\begin{equation*}
\mathcal{S}(s)_\alpha^0 = \bigcap_{\substack{{\bf v'} < {\bf v} \\ \slope_\alpha(Q, {\bf v'}) \geq \slope_\alpha(Q, {\bf v}) }} \mathcal{S}(s)_{\bf v'}^0
\end{equation*}
is residual in $\mathcal{S}(s)$. In other words, there is a residual set $\mathcal{S}(s)_\alpha^0$ of $C^1$ homotopy equivalences between $s$ and the constant map, such that for every $F \in \mathcal{S}(s)_\alpha^0$ the image $F( S^n \times [0,1])$ is contained in $\Rep(Q, {\bf v})^{\alpha-st}$.

Since residual sets are non-empty, then there is at least one basepoint-preserving homotopy equivalence in $\Rep(Q, {\bf v})^{\alpha-st}$ between $s$ and the constant map, and so $\pi_n (\Rep(Q, {\bf v})^{\alpha-st}) = 0$.
\end{proof}

\begin{corollary}\label{cor:homotopy-moduli-space}
If $n + 1 < d_{min} (Q, {\bf v}, \alpha)$ then $\pi_n \left(\mathcal{M}_\alpha(Q, {\bf v})^{st} \right) \cong \pi_{n-1} (PG_{\bf v})$.
\end{corollary}

\begin{proof}
The group $PG_{\bf v}$ acts freely on $\Rep(Q, {\bf v})^{\alpha-st}$, and the quotient is 
\begin{equation*}
\mathcal{M}_\alpha(Q, {\bf v})^{st} = \Rep(Q, {\bf v})^{\alpha-st} / PG_{\bf v} .
\end{equation*}
 The fibre bundle
\begin{equation*}
PG_{\bf v} \rightarrow \Rep(Q, {\bf v})^{\alpha-st} \rightarrow \mathcal{M}_\alpha (Q, {\bf v})^{st}
\end{equation*}
induces a long exact sequence of homotopy groups
\begin{multline*}
\cdots \rightarrow \pi_n \left(\Rep(Q, {\bf v})^{\alpha-st} \right) \rightarrow \pi_n \left(\mathcal{M}_\alpha(Q, {\bf v})^{st} \right) 
\rightarrow \pi_{n-1} (PG_{\bf v}) \rightarrow \pi_{n-1} \left( \Rep(Q, {\bf v})^{\alpha-st} \right) \rightarrow \cdots ,
\end{multline*}
and therefore an application of Theorem \ref{thm:homotopy-trivial} shows that 
$$
\pi_n \left(\Rep(Q, {\bf v})^{\alpha-st} \right) = \pi_{n-1} \left( \Rep(Q, {\bf v})^{\alpha-st} \right) = 0
$$
for $n < d_{min}(Q, {\bf v}, \alpha) - 1$. This shows that from the long exact sequence above, we have the isomorphism
\begin{equation*}
\pi_n \left(\mathcal{M}_\alpha(Q, {\bf v})^{st} \right) \cong \pi_{n-1} (PG_{\bf v}) ,
\end{equation*}
which completes the proof.
\end{proof}

\section{Examples}\label{sec:examples}

This section contains two applications of Theorem \ref{thm:homotopy-trivial} and Corollary \ref{cor:homotopy-moduli-space}. The first is to study the low-dimensional topology of moduli spaces that contain the stable loci of Nakajima quiver varieties, an explicit calculation is performed for the case of the ADHM quiver with the results in Theorem \ref{thm:doubled-adhm-homotopy}. The second application of the main theorem is to the moduli space of polygons in $\R^3$, which parametrises polygons in $\R^3$ with a fixed number of sides (determined by the quiver), and fixed sidelengths (determined by the choice of parameter), up to equivalence by rigid motions in $\R^3$. For almost all choices of parameter (we call these \emph{generic} parameters) there do not exist any polygons that lie in a line, and the moduli space is smooth. For non-generic parameters, \emph{degenerate polygons} (those that lie in a line) correspond to singularities in the moduli space $\mathcal{M}_\alpha(Q, {\bf v})$, since there are rotations of $\R^3$ that leave such polygons invariant. In Theorem \ref{thm:polygon-homotopy} we compute the low-dimensional homotopy groups of the space of \emph{non-degenerate polygons} (those that are not contained in a line), in the case where the parameter $\alpha$ is non-generic.

\subsection{The ``Doubled'' ADHM quiver}\label{subsec:ADHM}

It is well-known that the moduli space of framed instantons on $S^4$ has a description in terms of a finite-dimensional hyperk\"ahler quotient construction, known as the \emph{ADHM construction} (see \cite{ADHM78}, \cite{Donaldson84} for more details).  

This finite-dimensional hyperk\"ahler quotient construction is a special case of what is now known as a \emph{Nakajima quiver variety}, where a new quiver $\tilde{Q}$ is constructed by adding a ``dual'' edge to each edge on the original quiver $Q$, and the moduli space $\mathcal{M}_{(\alpha, \beta)}^{hk}(Q, {\bf v})$ is defined to be the hyperk\"ahler quotient of $\Rep(\tilde{Q}, {\bf v}) = T^* \Rep(Q, {\bf v})$ with parameters $\alpha$ and $\beta$ (see \cite{Nakajima94} for a complete description and \cite{Kronheimer89}, \cite{KronheimerNakajima90} for more examples). 

There is an inclusion of each of these hyperk\"ahler quotients $\mathcal{M}_{(\alpha, \beta)}^{hk}(Q, {\bf v}) \hookrightarrow \mathcal{M}_\alpha(\tilde{Q}, {\bf v})$ into the moduli space associated to the quiver $\tilde{Q}$, with the same dimension vector and stability parameter. The purpose of this section is to use Theorem \ref{thm:homotopy-trivial} to study the topology of the space $\mathcal{M}_0(\tilde{Q}, {\bf v})^{st}$, which contains $\mathcal{M}_{(0,0)}^{hk}(Q, {\bf v})$, in the case of the quiver corresponding to the ADHM construction of framed instantons on $S^4$. It is worth mentioning also that the method used is more general than just this one specific example, one can also prove analogous results for any Nakajima quiver variety by computing $d_{min}$ for each example.

Fix a vector bundle $E$ over $S^4$ with rank $n$ and instanton number $k$. Then the quiver $\tilde{Q}$ under consideration is given as follows.

\begin{equation}\label{eqn:ADHM-quiver}
\xygraph{
!{<0cm, 0cm>;<1.5cm, 0cm>:<0cm, 1.5cm>::}
!{(0,1.5) }*+{\bullet_{\C^k}}="a"
!{(0,0) }*+{\bullet_\C}="b"
"a" :@/^0.1cm/^\cdots "b" "a" :@/^0.5cm/ "b" "a" :@(l,lu) "a"
"b" :@/^0.1cm/^\cdots "a" "b" :@/^0.5cm/ "a" "a" :@(r,ru) "a"}
\end{equation}
(note that this differs slightly from Donaldson's description in \cite{Donaldson84}, here we use the construction of Crawley-Boevey in \cite{Crawley01} which identifies the moduli space associated to the quiver described above with the solutions to the equations from \cite{Donaldson84}).

Let $\alpha = 0$ be the stability parameter. There are $n$ arrows from the vertex $\C^k$ to the vertex $\C$, another $n$ arrows from the vertex $\C$ to the vertex $\C^k$, and two arrows from $\C^k$ to itself. In this case $G_{\bf v} = {\sf U}(k) \times {\sf U}(1)$, and $PG_{\bf v} \cong {\sf U}(k)$. 

\begin{remark}
When $\alpha = 0$ then every representation is $\alpha$-semistable, and a representation is $\alpha$-stable if and only if there are no non-trivial subrepresentations. 
\end{remark}

A subrepresentation $\pi$ can have one of the following forms (the diagram below also includes the quiver for $\pi^\perp$). Firstly we have the case where the vertex $\C$ is contained in the subrepresentation $\pi$. Since $\pi$ cannot have dimension vector ${\bf v}$, then $k_1 < k$ in the diagram below.
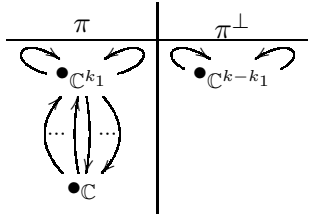
\begin{equation*}
\xygraph{
!{<0cm, 0cm>;<1cm, 0cm>:<0cm, 1cm>::}
!{(0,2.2) }*+{\pi}="p"
!{(2,2.2) }*+{\pi^\perp}="q"
!{(1,2.5) }*+{}="r"
!{(1,-0.3) }*+{}="s"
!{(3,2.0) }*+{}="u"
!{(-1,2.0) }*+{}="v"
!{(0,1.5) }*+{\bullet_{\C^{k_1}}}="a"
!{(0,0) }*+{\bullet_\C}="b"
!{(2,1.5) }*+{\bullet_{\C^{k-k_1}}}="c"
"a" :@/^0.1cm/^\cdots "b" "a" :@/^0.5cm/ "b" "a" :@(l,lu) "a"
"b" :@/^0.1cm/^\cdots "a" "b" :@/^0.5cm/ "a" "a" :@(r,ru) "a"
"c" :@(l,lu) "c" "c" :@(r,ru) "c" "r" - "s" "u" - "v"}
\end{equation*}
The elements of $\Hom^1(\pi, \pi^\perp)$ correspond to maps from the vertex $\C^{k_1}$ to $\C^{k-k_1}$ (there are two edges with this property) and maps from the vertex $\C$ to $\C^{k-k_1}$ (there are $n$ edges with this property). Therefore we see that $\dim_\R \Hom^1(\pi, \pi^\perp) = 4 k_1 (k - k_1) + 2 n (k - k_1)$. A similar computation shows that $\dim_\R \Hom^0(\pi, \pi^\perp) = 2k_1 (k - k_1)$. Therefore in this case
\begin{equation}
\dim_\R \Hom^1(\pi, \pi^\perp) - \dim_\R \Hom^0(\pi, \pi^\perp) = 2 (k_1 + n) (k - k_1) .
\end{equation}

The second case is when the vertex $\C$ is contained in $\pi^\perp$. Since $\pi$ cannot have dimension vector ${\bf 0}$, then $k_1 > 0$ in the diagram below.
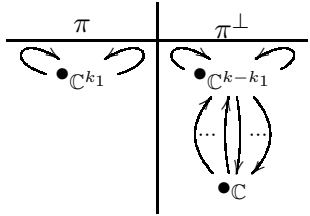
\begin{equation*}
\xygraph{
!{<0cm, 0cm>;<1cm, 0cm>:<0cm, 1cm>::}
!{(0,2.2) }*+{\pi}="p"
!{(2,2.2) }*+{\pi^\perp}="q"
!{(1,2.5) }*+{}="r"
!{(1,-0.3) }*+{}="s"
!{(3,2.0) }*+{}="u"
!{(-1,2.0) }*+{}="v"
!{(0,1.5) }*+{\bullet_{\C^{k_1}}}="a"
!{(2,0) }*+{\bullet_\C}="b"
!{(2,1.5) }*+{\bullet_{\C^{k-k_1}}}="c"
"c" :@/^0.1cm/^\cdots "b" "c" :@/^0.5cm/ "b" "a" :@(l,lu) "a"
"b" :@/^0.1cm/^\cdots "c" "b" :@/^0.5cm/ "c" "a" :@(r,ru) "a"
"c" :@(l,lu) "c" "c" :@(r,ru) "c" "r" - "s" "u" - "v"}
\end{equation*}
A calculation shows that 
$$
\dim_\R \Hom^1(\pi, \pi^\perp) = 4 k_1 (k - k_1) + 2 n k_1,
$$
and $\dim_\R \Hom^0(\pi, \pi^\perp) = 2k_1 (k - k_1)$. Therefore in this case
\begin{equation}
\dim_\R \Hom^1(\pi, \pi^\perp) - \dim_\R \Hom^0(\pi, \pi^\perp) = 2 k_1 (k - k_1 + n) .
\end{equation}

In both cases we see that $\dim_\R \Hom^1(\pi, \pi^\perp) - \dim_\R \Hom^0(\pi, \pi^\perp) \geq 2 (k-1 + n)$, and so $d_{min} \geq 2 (k-1+n)$. Therefore Theorem \ref{thm:homotopy-trivial} gives the following result.
\begin{theorem}\label{thm:doubled-adhm-homotopy}
Let $\tilde{Q}$ be the quiver in \eqref{eqn:ADHM-quiver}, and let $\alpha = 0$. Then
\begin{equation*}
\pi_j \left( \Rep(\tilde{Q}, {\bf v})^{\alpha-st} \right) = 0
\end{equation*}
for all  $j < 2 (k-1 + n) - 1$, and 
\begin{equation*}
\pi_j \left( \mathcal{M}_\alpha ( \tilde{Q}, {\bf v})^{st} \right) = \pi_j \left({\sf U}(k) \right)
\end{equation*}
for all $j < 2 (k-1 + n) -1 $.
\end{theorem}

\subsection{The moduli space of non-degenerate polygons in $\R^3$}\label{subsec:polygon}

The topology and geometry of the moduli space of polygons in $\R^3$ with fixed sidelengths has been well-studied, see for example \cite{Klyachko94}, \cite{KapovichMillson96}, \cite{HausmannKnutson97}, \cite{HausmannKnutson98}. When the sidelengths are chosen so that some polygons may lie in a line (we call these \emph{degenerate} polygons) then the moduli space has singularities, corresponding to the fact that a degenerate polygon is invariant under a circle subgroup of the rotation group. The purpose of this section is to prove new results about the low-dimensional topology of the moduli space of \emph{non-degenerate} polygons.

An $n$-sided polygon in $\R^3$ corresponds to a choice of basepoint $b \in \R^3$, as well as $n$ vectors $e_1, \ldots, e_n$ such that $e_1 + \cdots + e_n = 0$, where the vertices of the polygon are then given by $v_k = b + e_1 + \cdots + e_k$. The \emph{moduli space of polygons in $\R^3$ with fixed sidelengths} is the space of all polygons in $\R^3$ with basepoint $0$ and given sidelengths $s_1, \ldots, s_n$, modulo rotations by $\SO(3)$. It is well-known (see for example \cite{Konno02}) that this moduli space corresponds to the representation variety of a quiver, given as follows. Let $Q$ be the ``star-shaped'' quiver below, with one central vertex and $n$ outer vertices (one for each side of the polygon). Fix a dimension vector ${\bf v}$, with dimension $2$ at the central vertex, and dimension $1$ at each of the outer vertices. In this case $G_{\bf v} = {\sf U}(2) \times {\sf U}(1)^n$, and $PG_{\bf v} \cong {\sf U}(2) \times {\sf U}(1)^{n-1}$.

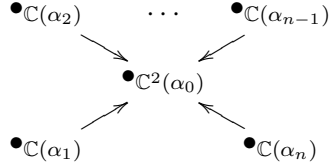
\begin{equation*}
\xygraph{
!{<0cm, 0cm>;<0.9cm, 0cm>:<0cm, 0.9cm>::}
!{(-1.73,1) }*+{\bullet_{\C (\alpha_2)}}="a"
!{(0,1) }*+{\cdots}="b"
!{(1.73,1) }*+{\bullet_{\C (\alpha_{n-1})}}="c"
!{(1.73,-1) }*+{\bullet_{\C (\alpha_n)}}="d"
!{(-1.73,-1) }*+{\bullet_{\C (\alpha_1)}}="e"
!{(0,0) }*+{\bullet_{\C^2 (\alpha_0)}}="z"
"a" : "z"   "c" : "z" "d" : "z" "e" : "z"}
\end{equation*}

With the conventions of Section \ref{sec:background}, the parameter $\alpha_j$ at the $j^{th}$ vertex is chosen so that $\alpha_j = -s_j$ (where $s_j$ is the sidelength of the $j^{th}$ edge of the polygon in $\R^3$), and the parameter $\alpha_0$ is then given by $\displaystyle{ \alpha_0 = -\frac{1}{2}\sum_{j=1}^n \alpha_j }$, i.e. the associated central element $\alpha \in \mathfrak{g}_{\bf v}^*$ is trace-free.

Since $\alpha_j < 0$ for each $j=1, \ldots, n$ then by the definition of slope stability in \eqref{eqn:def-slope-stability} we see that a representation which is not $\alpha$-stable must correspond to one of the following two cases. The first case consists of those subrepresentations $\pi$ that do not contain any subspaces of the central vertex (as for the examples in the previous section, the diagram below includes the subrepresentation $\pi^\perp$). The outer vertices in the following diagram have been re-labelled, so that $\alpha_{j_i}$ for $1 \leq i \leq k$ are in the subrepresentation $\pi$, and $\alpha_{j_i}$ for $k+1 \leq i \leq n$ are in the subrepresentation $\pi^\perp$.
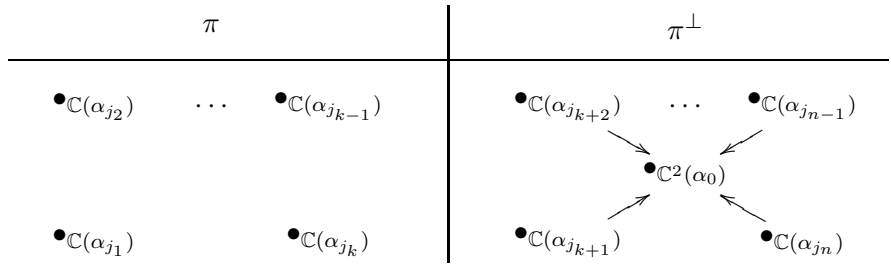
\begin{equation*}
\xygraph{
!{<0cm, 0cm>;<0.9cm, 0cm>:<0cm, 0.9cm>::}
!{(0,2.2) }*+{\pi}="p"
!{(7,2.2) }*+{\pi^\perp}="q"
!{(3.5,2.5) }*+{}="r"
!{(3.5,-1.3) }*+{}="s"
!{(10,1.7) }*+{}="u"
!{(-3,1.7) }*+{}="v"
!{(-1.73,1) }*+{\bullet_{\C (\alpha_{j_2})}}="a"
!{(0,1) }*+{\cdots}="b"
!{(1.73,1) }*+{\bullet_{\C (\alpha_{j_{k-1}})}}="c"
!{(1.73,-1) }*+{\bullet_{\C (\alpha_{j_k})}}="d"
!{(-1.73,-1) }*+{\bullet_{\C (\alpha_{j_1})}}="e"
!{(5.27,1) }*+{\bullet_{\C (\alpha_{j_{k+2}})}}="g"
!{(7,1) }*+{\cdots}="h"
!{(8.73,1) }*+{\bullet_{\C (\alpha_{j_{n-1}})}}="i"
!{(8.73,-1) }*+{\bullet_{\C (\alpha_{j_n})}}="j"
!{(5.27,-1) }*+{\bullet_{\C (\alpha_{j_{k+1}})}}="k"
!{(7,0) }*+{\bullet_{\C^2 (\alpha_0)}}="z"
"g" : "z" "i" : "z" "j" : "z" "k" : "z" "r" - "s" "u" - "v"}
\end{equation*}
Note that since the subrepresentation $\pi$ is non-trivial, then $k \geq 1$. The elements of $\Hom^1(\pi, \pi^\perp)$ correspond to maps from the vertices in the quiver on the left (corresponding to the subrepresentation $\pi$) to the central vertex in the quiver on the right (corresponding to $\pi^\perp$). A calculation shows that $\dim_\R \Hom^1 (\pi, \pi^\perp) = 4k$, and that $\dim_\R \Hom^0(\pi, \pi^\perp) = 0$. Therefore, in this case
\begin{equation}\label{eqn:non-split-centre}
\dim_\R \Hom^1 (\pi, \pi^\perp) - \dim_\R \Hom^0(\pi, \pi^\perp) = 4k \geq 4.
\end{equation}

The other possible type of subrepresentation occurs when the subrepresentation $\pi$ contains a complex dimension one subspace of the central vertex, described in the diagram below.

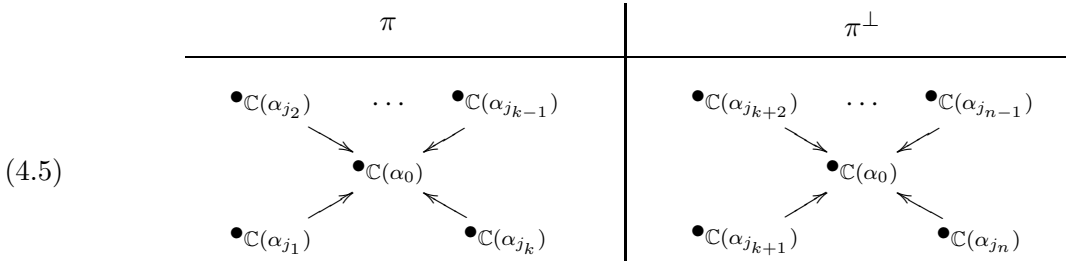
\begin{equation}\label{eqn:polygon-subrep}
\xygraph{
!{<0cm, 0cm>;<0.9cm, 0cm>:<0cm, 0.9cm>::}
!{(0,2.2) }*+{\pi}="p"
!{(7,2.2) }*+{\pi^\perp}="q"
!{(3.5,2.5) }*+{}="r"
!{(3.5,-1.3) }*+{}="s"
!{(10,1.7) }*+{}="u"
!{(-3,1.7) }*+{}="v"
!{(-1.73,1) }*+{\bullet_{\C (\alpha_{j_2})}}="a"
!{(0,1) }*+{\cdots}="b"
!{(1.73,1) }*+{\bullet_{\C (\alpha_{j_{k-1}})}}="c"
!{(1.73,-1) }*+{\bullet_{\C (\alpha_{j_k})}}="d"
!{(-1.73,-1) }*+{\bullet_{\C (\alpha_{j_1})}}="e"
!{(0,0) }*+{\bullet_{\C (\alpha_0)}}="f"
!{(5.27,1) }*+{\bullet_{\C (\alpha_{j_{k+2}})}}="g"
!{(7,1) }*+{\cdots}="h"
!{(8.73,1) }*+{\bullet_{\C (\alpha_{j_{n-1}})}}="i"
!{(8.73,-1) }*+{\bullet_{\C (\alpha_{j_n})}}="j"
!{(5.27,-1) }*+{\bullet_{\C (\alpha_{j_{k+1}})}}="k"
!{(7,0) }*+{\bullet_{\C (\alpha_0)}}="z"
"a" : "f" "c" : "f" "d" : "f" "e" : "f" 
"g" : "z" "i" : "z" "j" : "z" "k" : "z" "r" - "s" "u" - "v"}
\end{equation}
Unlike the previous case (where any value of $k$ between $1$ and $n$ could occur), the range of possible values of $k$ is limited by the choice of $\alpha_j$ at each vertex (i.e. the sidelengths of the polygon). To see this, let ${\bf v'}$ be the dimension vector associated to the subrepresentation $\pi$, and note that for the condition $\slope_\alpha(Q, {\bf v'}) \geq \slope_\alpha(Q, {\bf v})$ to hold we must have
\begin{equation*}
\deg(Q, {\bf v'}) = -\alpha_0 - \sum_{i=1}^k \alpha_{j_i} \geq 0 ,
\end{equation*}
which places a restriction on the possible configurations that can occur in this case (recall that $\alpha_0 > 0$ and that $\alpha_j < 0$ for $j = 1, \ldots, k$.

A calculation shows that $\dim_\R \Hom^1(\pi, \pi^\perp) = 2k$, and that $\dim_\R \Hom^0(\pi, \pi^\perp) = 2$. Therefore, in this case
\begin{equation}\label{eqn:split-centre}
\dim_\R \Hom^1 (\pi, \pi^\perp) - \dim_\R \Hom^0(\pi, \pi^\perp) = 2k-2 .
\end{equation}

In analogy with \cite{HausmannKnutson98}, define a subset $J \subset \{ 1, \ldots, n \}$ to be \emph{strictly short} if 
\begin{equation}
-\sum_{j \in J} \alpha_j < -\sum_{j \notin J} \alpha_j .
\end{equation}
Equivalently, $J$ is strictly short if $-\alpha_0 - \sum_{j \in J} \alpha_j < 0$.

\begin{definition}
For a positive integer $\ell$, the parameter $\alpha$ is \emph{$\ell$-strictly short} if every subset $J \subset \{ 1, \ldots, n \}$ with $\ell$ elements is strictly short. If there is no positive integer $\ell$ for which $\alpha$ is $\ell$-strictly short then $\alpha$ is called \emph{$0$-strictly short}.
\end{definition}

\begin{remark}
If $\alpha$ is $\ell$-strictly short, then there are no subrepresentations of the type in diagram \eqref{eqn:polygon-subrep} with $k \leq \ell$. Therefore by \eqref{eqn:split-centre} we have
\begin{equation}\label{eqn:ell-short-bound}
\dim_\R \Hom^1 (\pi, \pi^\perp) - \dim_\R \Hom^0(\pi, \pi^\perp) \geq 2\ell .
\end{equation}
\end{remark}

Equations \eqref{eqn:non-split-centre} and \eqref{eqn:ell-short-bound} show that $d_{min}$ takes the following values when $\alpha$ is $\ell$-strictly short.
\begin{center}
\begin{tabular}{l | c} 
$d_{min}$ & $\ell$ \\ \hline
0 & 0 \\
2 & 1 \\
4 & 2 or greater
\end{tabular}
\end{center}
Together with Theorem \ref{thm:homotopy-trivial} and Corollary \ref{cor:homotopy-moduli-space}, the above calculations lead to the following theorem.
\begin{theorem}\label{thm:polygon-homotopy}
\begin{enumerate}

\item If $\alpha$ is $2$-strictly short then 

\begin{itemize}
\item $\pi_j \left(\Rep(Q, {\bf v})^{\alpha-st} \right) = 0$ for $j \leq 2$.

\item $\pi_j \left( \mathcal{M}_\alpha (Q, {\bf v})^{st} \right) = \pi_j \left( \U(2) \times \U(1)^{n-1} \right)$ for $j \leq 2$.
\end{itemize}

\item If $\alpha$ is $1$-strictly short then $\mathcal{M}_\alpha (Q, {\bf v})^{st}$ and $\Rep(Q, {\bf v})^{\alpha-st}$ are both connected.
\end{enumerate}
\end{theorem}

\section*{Acknowledgments}
The author would like to thank Georgios Daskalopoulos for his advice and encouragement, as well as the referee for some useful observations and background references.

\def\cprime{$'$}

\end{document}